\pdfoutput=1
\documentclass[preprint,12pt]{elsarticle}

\usepackage{amssymb}
\usepackage{amsmath}
\usepackage{graphicx}
\usepackage{booktabs}
\usepackage{algorithm}
\usepackage{algpseudocode}
\usepackage{hyperref}
\usepackage{geometry}

\journal{arXiv}

\begin{document}

\begin{frontmatter}

\title{A Real-Time Scalable Heuristic DSS Framework for Capacity-Constrained Retail Allocation under Supply Chain Uncertainty}

\author[1]{Abd\"ussamet S\"okel\corref{cor1}}
\ead{samet.soekel@yandex.com}
\cortext[cor1]{Corresponding author. Postal address: Industrial Engineering Department, Eskisehir Osmangazi University, Turkey.}

\address[1]{Industrial Engineering Department, Eskisehir Osmangazi University, Turkey}

\begin{abstract}
The rapid proliferation of omnichannel retail strategies has fundamentally transformed store replenishment operations in uncertain supply chain environments. With retail stores increasingly acting as hybrid fulfillment centers, pooled inventory allocation must absorb uncertain order realizations, constrained receiving capacities, dynamic vehicle limits, multi-tiered product priorities, and planner-controlled outbound warehouse preferences. This study frames this commercial reality as an extended constrained variant of the Multidimensional Knapsack Problem (MKP). Recognizing that exact optimization techniques such as Mixed-Integer Linear Programming (MILP) are computationally prohibitive in large-scale real-time settings, we propose a real-time scalable heuristic embedded in a computationally efficient Decision Support System (DSS) framework based on set-oriented cumulative filtering. The framework evaluates cumulative flow-through deductions, third-party logistics routing integrations, category-specific volume caps, warehouse activation filters, and user-defined warehouse priority ranks. An extensive case study within a large retail network covering 212{,}278 order records from June 2025 to April 2026 demonstrates the impact of the proposed methodology. Using January 2026 as the go-live cutoff, weighted ship-to-order ratio improved from 54.1\% to 67.8\%, weighted same-day coverage improved from 24.3\% to 37.8\%, and store-days with order volumes above store limits were reduced by 48.6\%. These findings indicate that the proposed real-time scalable heuristic and computationally efficient DSS framework provide practical, uncertainty-aware allocation support for volatile retail supply chains.
\end{abstract}


\end{frontmatter}

\section{Introduction}
\label{sec:intro}

The retail industry continues to undergo a paradigm shift driven by evolving supply chain architectures and fluctuating consumer demands. Traditional continuous-review replenishment models, which historically relied on stable lead times and infinite store receiving capacities, are increasingly insufficient \cite{Boysen2021}. In the contemporary omnichannel landscape, "ship-from-store" and "buy-online-pickup-in-store" (BOPIS) fulfillment strategies necessitate that brick-and-mortar stores function not only as points of sale but as active micro-fulfillment centers \cite{Hubner2022}. Consequently, dynamic receiving capacities heavily restricted by inbound returns, cross-docking operations, and backroom storage constraints have emerged as a critical bottleneck.

Efficient inventory allocation from a central Distribution Center (DC) to a vast retail network is a highly constrained resource allocation challenge. Retailers must rapidly assign thousands of discrete product variants to heterogeneous stores under strict daily payload thresholds. In practice, this allocation frequently starts from multiple outbound warehouses, and planners may need to activate only a subset of warehouses and enforce a strict precedence order among them for a given planning cycle. Additionally, logistics practitioners face nested constraints: specific product categories (e.g., bulky or fragile items) possess separate sub-quotas per transportation route, while routing logic must distinguish between internal fleets and outsourced third-party logistics (3PL) providers based on real-time eligibility matrices. These constraints are further intensified by uncertainty in requested quantities, residual store capacities, and route-level feasible payload windows.

Analytically, this problem maps to the Generalized Assignment Problem (GAP) and the Multidimensional Knapsack Problem (MKP), both known to be strongly NP-hard \cite{Kellerer2004}. While academic literature extensively applies exact solvers via Mixed-Integer Linear Programming (MILP), they suffer from severe intractability in commercial environments. In recent omnichannel and uncertainty-oriented formulations, this computational burden remains a central barrier and generally requires decomposition or heuristic acceleration layers for operational-scale decision cycles \cite{Bansal2024,Qiu2025,Roushdy2026,Bahuguna2026}.

To bridge this operational void, we propose a highly efficient, data-driven real-time scalable heuristic tailored for decision support systems. Unlike iterative meta-heuristics that exhibit lengthy convergence times, the proposed methodology leverages set-based operations and relational algebra to enforce knapsack thresholds instantly. The core contribution of this research is a computationally efficient DSS framework that avoids slow sequential iterative algorithms while introducing operator-level warehouse selection and prioritization controls.

This paper contributes along five dimensions: it formalizes a multi-constraint allocation problem with operational uncertainty, develops a real-time scalable heuristic for day-level deployment, derives analytical properties of the heuristic (feasibility preservation and $\mathcal{O}(N\log N)$ computational complexity), embeds the procedure in a computationally efficient DSS framework for industrial execution, and reports empirical evidence from a full-scale retailer before-after implementation.

\section{Literature Review}
\label{sec:litreview}

To address the complexities of modern retail allocation, our research draws upon five tightly connected domains of literature: omnichannel replenishment and fulfillment coupling, omnichannel logistics network design, scalable hybrid exact-heuristic optimization, uncertainty-oriented decision frameworks in supply chains, and modern extensions of MKP-based constrained allocation models.

\subsection{Omnichannel Retail Replenishment and Dynamic Constraints}
Recent Q1/Q2 operations and logistics literature shows that omnichannel retailing transforms replenishment into a coupled, dynamic decision problem rather than a static store-capacity check. Survey evidence confirms this structural shift and the operational pressure on stores acting as hybrid fulfillment nodes \cite{Boysen2021,Hubner2022}. At policy level, dynamic ship-from-store allocation and in-store fulfillment process design have been formalized under demand and cost uncertainty \cite{Bayram2021,Difrancesco2021}. More recent studies explicitly integrate replenishment and online-demand allocation across periods, showing that myopic rules can create inventory imbalance and service degradation \cite{Goedhart2023,Bansal2024}.

In parallel, strategic network-design studies extend the operational view by embedding multi-tier, multi-service-level, and multi-modal delivery architecture into omnichannel planning \cite{Arslan2021,Janjevic2021,Janjevic2020}. Customer-facing delivery/return channel preferences and walking-distance constraints are also modeled as endogenous design drivers in last-mile structures \cite{GuerreroLorente2020}. Together, these findings support a dynamic and network-aware perspective where fulfillment responsiveness, not only inventory depth, determines feasible allocation quality.

\subsection{Heuristics in Real-Time Logistics}
Computational realism is a recurring theme in recent omnichannel optimization research. Although exact or stochastic formulations deliver strong policy quality, practical deployment typically depends on decomposition and fast heuristics layered on top of mathematical programs. For example, integrated batching-routing operations in grocery contexts are solved via adaptive large neighborhood search to preserve tractability at realistic scales \cite{Kuhn2021}. Scenario-reduced two-stage stochastic models combined with dynamic fulfillment heuristics have also shown significant cost improvements in uncertain omnichannel settings \cite{Abouelrous2022}. Likewise, demand-uncertain replenishment with service-level constraints has been treated through rule-structured MILP frameworks and multi-layer solution procedures to maintain decision speed \cite{Qiu2025}. Joint pricing-fulfillment models in multi-zone networks similarly rely on linearization to obtain implementable mixed-integer programs \cite{Pichka2022}.

This stream converges on a practical insight: robust industrial decision support tends to favor hybrid methods that preserve core optimization structure while enforcing deterministic and short runtime behavior.

\subsection{Uncertainty-Oriented Decision Frameworks in Supply Chains}
Recent uncertainty-focused SCM studies demonstrate growing interest in uncertainty-aware model design, especially for replenishment coordination, reverse logistics, and sensitivity-driven managerial analysis \cite{Roushdy2026,Bahuguna2026}. These studies highlight that incorporating uncertainty is not only a modeling preference but an operational necessity when demand regimes, return quality, and route feasibility change across planning cycles. Our work complements this stream by emphasizing a deployment-oriented perspective: rather than introducing a high-dimensional stochastic solver, we design a real-time scalable heuristic and computationally efficient DSS framework that can be rerun quickly under refreshed daily realizations.

\subsection{Multidimensional Knapsack Problem in Retail}
The constrained allocation core of the problem remains grounded in classical MKP/GAP theory \cite{Kellerer2004}. However, contemporary omnichannel formulations add operational dimensions not present in canonical knapsack settings, including integrated replenishment-allocation coupling \cite{Bansal2024}, channel-specific service-level control under uncertainty \cite{Qiu2025}, and responsiveness-driven decentralization of inventory deployment \cite{Kim2023Network}. Despite this progress, a gap remains in models that simultaneously represent nested route-category quotas, real-time store residual capacities, route eligibility matrices, and planner-controlled warehouse precedence in one executable decision pipeline. Our formulation targets this gap by combining these constraints in a deterministic cumulative-allocation framework suitable for real-time DSS execution.

\section{Problem Definition and Mathematical Formulation}
\label{sec:mathmodel}

The challenge is to select a subset of pending orders from a central pool destined for physical stores, maximizing prioritized dispatched value while ensuring that no multi-tier capacities are breached and that planner-defined outbound warehouse selection and precedence are respected.

\subsection{Indices and Parameters}
\begin{itemize}
    \item $M$: Set of all retail stores, indexed by $m$.
    \item $R$: Set of delivery route aggregates, indexed by $r$.
    \item $P$: Set of product categories, indexed by $p$.
    \item $W$: Set of outbound warehouses, indexed by $w$.
    \item $I$: Set of unallocated pending orders, indexed by $i$.
    \item $v_i$: Dimensional weight or discrete volume (e.g., Desi) of order $i$.
    \item $w_i$: Business priority weight of order $i$.
    \item $C_m^{(t)}$: Dynamic residual receiving capacity of store $m$ at planning day $t$, computed from base capacity after flow-through and accepted-load deductions.
    \item $L_p$: Cumulative capacity limit associated with sensitive product category $p$.
    \item $E_{m,p} \in \{0,1\}$: Binary eligibility matrix indicating whether store $m$ can receive product group $p$.
    \item $z_w \in \{0,1\}$: Warehouse activation parameter (1 if warehouse $w$ is selected for the run).
    \item $\pi_w \in \{1,\dots,\Pi_{\max}\}$: User-defined warehouse priority rank (smaller is higher priority).
\end{itemize}

\subsection{Decision Variables}
\begin{itemize}
    \item $X_i \in \{0,1\}$: 1 if order $i$ is allocated, 0 otherwise.
\end{itemize}

\subsection{MILP Formulation}

The objective function maximizes prioritized throughput while encoding warehouse precedence as a dominant term:
\begin{equation}
\max \sum_{i \in I} \left(\lambda(\Pi_{\max}+1-\pi_{w(i)}) + w_i\right)X_i
\end{equation}
where $\lambda$ is selected sufficiently large to preserve lexicographic preference for warehouse rank before business priority.

Subject to the following operational constraints:

1. \textbf{Dynamic Store Capacity:}
\begin{equation}
\sum_{i \in I: m(i)=m} v_i X_i \leq C_m^{(t)} \quad \forall m \in M
\end{equation}

2. \textbf{Product Category Limits per Route:} Let $R(m)$ represent the assigned route of store $m$, and $p(i)$ represent the category of order $i$.
\begin{equation}
\sum_{i \in I: p(i)=p, R(m(i))=r} v_i X_i \leq L_p \quad \forall p \in P_{\text{constrained}}, \forall r \in R
\end{equation}

3. \textbf{Node Eligibility Filtering:}
\begin{equation}
X_i \leq E_{m(i),p(i)} \quad \forall i \in I
\end{equation}

4. \textbf{Outbound Warehouse Activation:}
\begin{equation}
X_i \leq z_{w(i)} \quad \forall i \in I
\end{equation}

5. \textbf{Binary Constraint:}
\begin{equation}
X_i \in \{0,1\} \quad \forall i \in I
\end{equation}

\section{The Proposed Methodology: A Real-Time Scalable Heuristic}
\label{sec:heuristic}

Mathematical solvers attempt to evaluate the entire combinatorial space to resolve equations (1)-(6). The proposed real-time scalable heuristic operates differently by leveraging warehouse-aware lexicographic ranking and cumulative aggregate filtering, reducing complexity to $\mathcal{O}(N \log N)$ where $N$ is the number of pending orders.

The algorithm structures data into a sequential pipeline. First, invalid pathways are pruned by node eligibility and warehouse activation ($E_{m(i),p(i)} = 0$ or $z_{w(i)} = 0$). Next, remaining orders are sorted by warehouse rank first ($\pi_{w(i)}$), then by business priority. Using set-based processing, the system calculates continuous cumulative totals over warehouse, store, and route partitions. In the second-day simulation pass, store capacities are recalculated with additional accepted-load deductions before rerunning the same ordered filtering logic.

\subsection{Operational Uncertainty Representation}
The deployment setting treats key inputs as realized uncertain variables at each planning run: the order pool composition, dynamic residual capacities $C_m^{(t)}$, route-category exposure, and eligibility status. Instead of solving a static long-horizon model once, the DSS refreshes these inputs and reruns the real-time scalable heuristic in rolling cycles. This design provides a practical uncertainty-handling mechanism that is computationally tractable for daily operational use.

\subsection{Theoretical Properties of the Heuristic}
\noindent\textbf{Proposition 1 (Feasibility preservation).} Every allocation returned by the heuristic satisfies constraints (2)--(6) of the formulation.

\noindent\textit{Proof sketch.} Orders violating eligibility or warehouse activation are removed before ranking, which enforces (4)--(5). During sequential acceptance, each candidate is admitted only if store-level cumulative load remains within $C_m^{(t)}$, enforcing (2). For constrained categories, route-category cumulative load is checked against $L_p$ before acceptance, enforcing (3). Since accepted orders are binary selections from eligible candidates, (6) holds directly. Therefore, the final accepted set is feasible with respect to the modeled constraints.

\noindent\textbf{Proposition 2 (Computational complexity).} Let $N$ denote the number of eligible orders in a planning cycle. The heuristic runs in $\mathcal{O}(N\log N)$ time.

\noindent\textit{Proof sketch.} Eligibility filtering is linear in input size. The dominant step is sorting by warehouse rank and priority, costing $\mathcal{O}(N\log N)$. The subsequent cumulative pass over sorted orders updates constant-time load trackers and therefore costs $\mathcal{O}(N)$. The total complexity is thus $\mathcal{O}(N\log N)$.

\begin{algorithm}
\caption{Warehouse-Aware Cumulative Allocation Heuristic}
\begin{algorithmic}[1]
\State \textbf{Input:} Set $I, C_m^{(t)}, L_p, E, z, \pi$
\State \textbf{Output:} Allocated set $A$
\State $A \gets \emptyset$
\State $I_{eligible} \gets \{i \in I \mid E_{m(i),p(i)} = 1 \land z_{w(i)} = 1\}$
\State Sort $I_{eligible}$ by $\pi_{w(i)}$ ascending, then $w_i$ descending, then $v_i$ descending
\For{each store $m \in M$}
    \State $\text{StoreLoad}[m] \gets 0$
\EndFor
\For{each route $r \in R$ and category $p$}
    \State $\text{CategoryLoad}[r,p] \gets 0$
\EndFor
\For{each order $i \in I_{eligible}$ sequentially}
    \State $m \gets m(i)$, $p \gets p(i)$, $r \gets R(m)$
    \If{$\text{StoreLoad}[m] + v_i \leq C_m^{(t)}$}
        \If{category $p$ is constrained \textbf{and} $\text{CategoryLoad}[r,p] + v_i > L_p$}
            \State \textbf{continue}
        \Else
            \State Add $i$ to $A$
            \State $\text{StoreLoad}[m] \gets \text{StoreLoad}[m] + v_i$
            \If{category $p$ is constrained}
                \State $\text{CategoryLoad}[r,p] \gets \text{CategoryLoad}[r,p] + v_i$
            \EndIf
        \EndIf
    \EndIf
\EndFor
\State \textbf{return} $A$
\end{algorithmic}
\end{algorithm}

\subsection{Step-by-Step Execution Logic}
For reproducibility, the same procedure can be summarized in deterministic operational steps:
\begin{enumerate}
    \item Build the eligible order set by removing records that violate node eligibility or warehouse activation rules.
    \item Sort eligible orders by warehouse priority rank, then business priority, then order volume.
    \item Initialize store-level cumulative loads to zero.
    \item Initialize route-category cumulative loads to zero for constrained categories.
    \item Iterate through sorted orders sequentially.
    \item For each candidate, first check whether adding the order violates store residual capacity.
    \item If the order belongs to a constrained category, check route-category cumulative limit feasibility.
    \item Accept only orders that pass all active checks, and immediately update cumulative load trackers.
    \item Return the accepted order set as the final dispatch plan for the run.
\end{enumerate}

\section{Computationally Efficient DSS Framework and Real-World Implementation}
\label{sec:casestudy}

The proposed framework was deployed within the distribution architecture of a large omnichannel retailer. The operational environment couples an interactive Decision Support System with a centralized enterprise database, enabling large-scale logical aggregates over relational structures.

The interface exposes two operational controls before each run: (i) selectable outbound warehouses and (ii) drag-and-drop warehouse ordering to define dispatch precedence. The selected warehouses and ranks are materialized into a temporary priority table and inner-joined with the order pool, so non-selected warehouses are excluded at source.

The technical core bypasses traditional row-by-row loops. Constraint evaluation is executed through window-based cumulative data operations in the relational layer, with ordering logic starting from warehouse priority and followed by business hierarchy keys. For the second-day simulation, capacities are recalculated with accepted-load deductions and reused in the same ranking pipeline. Downstream export logic produces warehouse-specific output files, preserving traceability from decision stage to execution stage. In the results section, warehouses are reported using role labels (Warehouse-Primary, Warehouse-Auxiliary, Warehouse-InnerProducts), carriers as Carrier-XX, and stores as Store-XXXX to keep presentation concise and consistent.

\section{Experimental Results}
\label{sec:results}

We evaluated the decision support system with a before-after design using January 2026 as the go-live cutoff. The input data includes 212{,}278 order lines spanning June 2025 to April 2026, with three outbound warehouses and 772 destination stores. Throughout this section, \textit{batch} denotes the operational shipment unit. To avoid metric inflation in days where shipped batches can exceed same-day requested batches, we report two complementary service metrics:
\begin{itemize}
    \item \textbf{Ship/Order Ratio:} $\sum S / \sum R$ where $S$ is shipped batch and $R$ is requested batch (captures total outbound throughput).
    \item \textbf{Same-Day Coverage:} $\sum \min(S,R) / \sum R$ (bounded demand coverage proxy).
\end{itemize}

\begin{table}[htbp]
\centering
\caption{Before-After KPI Comparison (Go-Live: January 2026)}
\label{tab:results}
\begin{tabular}{@{}lccc@{}}
\toprule
\textbf{Metric} & \textbf{Before} & \textbf{After} & \textbf{Delta} \\ \midrule
Analysis horizon (days) & 188 & 74 & -- \\
Weighted Ship/Order Ratio & 54.1\% & 67.8\% & +13.7 pp \\
Weighted Same-Day Coverage & 24.3\% & 37.8\% & +13.4 pp \\
Average Daily Unserved Batch & 6,657 & 5,137 & -22.8\% \\
Store-days with Order $>$ Limit & 4.53\% & 2.33\% & -48.6\% \\
Store-days with Ship $>$ Limit & 0.82\% & 0.49\% & -40.6\% \\
Store-days with Full Fulfillment & 14.1\% & 20.9\% & +6.7 pp \\
\bottomrule
\end{tabular}
\end{table}

The increase in daily weighted same-day coverage is statistically significant (Mann--Whitney U, $p=6.0\times10^{-7}$), indicating that post go-live gains are unlikely to be explained by random day-level variation.

\noindent\textbf{Baseline choice.} The empirical benchmark is the pre go-live operational policy, which is the most relevant managerial baseline for this deployment context. This comparison directly measures decision-quality changes attributable to the proposed real-time scalable heuristic and computationally efficient DSS framework under the same operational environment.

\subsection{Method-Level Comparative Positioning}
In addition to the before-after benchmark, Table~\ref{tab:method-positioning} positions the proposed framework against two common alternatives in this literature: exact optimization-oriented pipelines and meta-heuristic search pipelines \cite{Qiu2025,Kuhn2021}. This is an analytical positioning table rather than a head-to-head computational experiment due to real-time operational constraints and industrial data restrictions.

\begin{table}[htbp]
\centering
\caption{Analytical positioning versus common solution families}
\label{tab:method-positioning}
\begin{tabular}{@{}p{3.6cm}p{3.4cm}p{3.4cm}p{4.0cm}@{}}
\toprule
\textbf{Criterion} & \textbf{Exact/Math-program-heavy} & \textbf{Meta-heuristic-heavy} & \textbf{Proposed framework} \\ \midrule
Daily rerun latency behavior & Often sensitive to instance size and solver tuning & Better than exact in many cases but still iteration-dependent & Deterministic ordered filtering with short runtime profile \\
Nested operational constraints (store, route, category, eligibility) & High modeling flexibility but computational burden can be high & Usually feasible with tailored encoding and penalty calibration & Explicitly embedded as hard cumulative checks in one pass \\
Planner controls (warehouse activation + rank) & Possible but often requires model refactoring between runs & Usually custom-coded and parameter-dependent & Native UI-driven controls directly mapped to ranking and filters \\
Operational auditability & Moderate (solver logs, model complexity) & Moderate to low (search trajectory dependence) & High (rule-transparent sequence and cumulative thresholds) \\
\bottomrule
\end{tabular}
\end{table}

\subsection{Temporal Dynamics and Throughput}
\begin{figure}[htbp]
\centering
\includegraphics[width=0.98\linewidth]{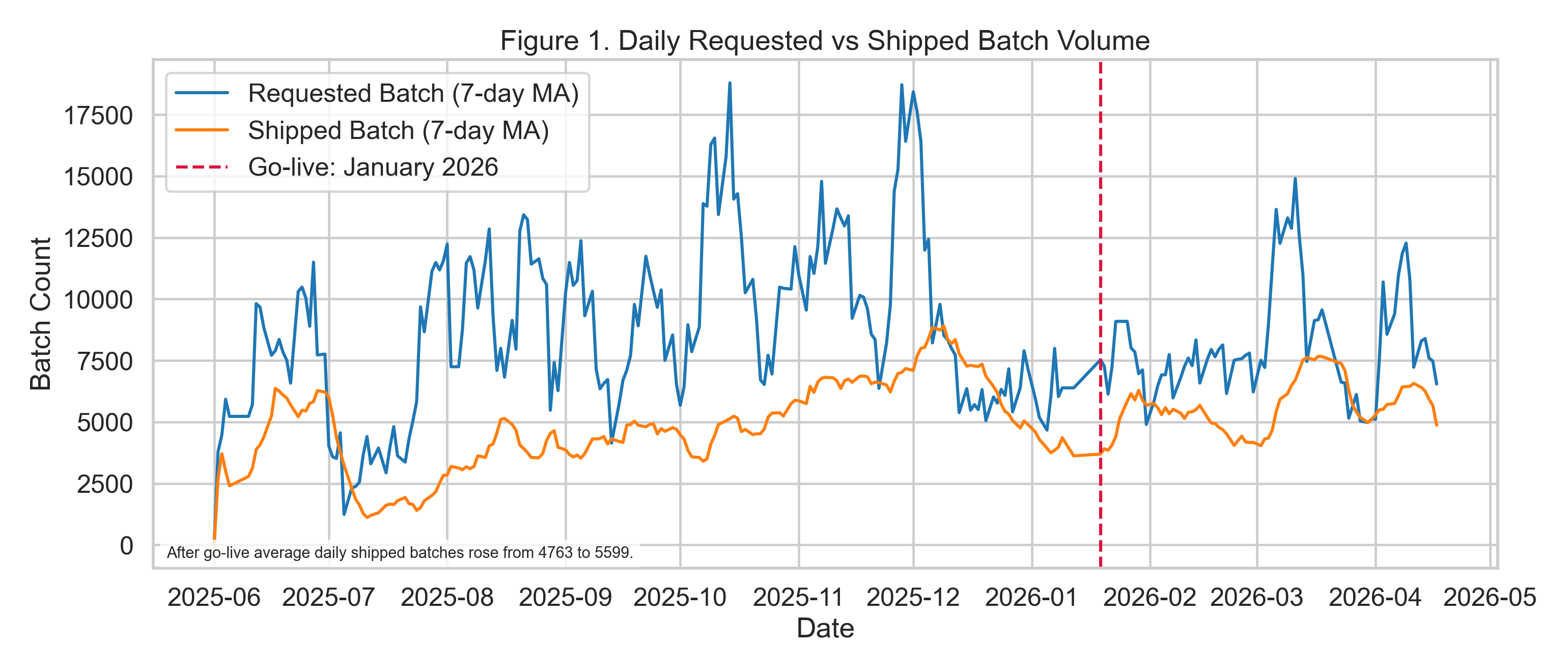}
\caption{Daily requested and shipped batch trajectories (7-day moving averages). The dashed vertical line marks go-live (January 2026). Post go-live, shipped volume follows requested volume more closely.}
\label{fig:daily-volumes}
\end{figure}

\begin{figure}[htbp]
\centering
\includegraphics[width=0.98\linewidth]{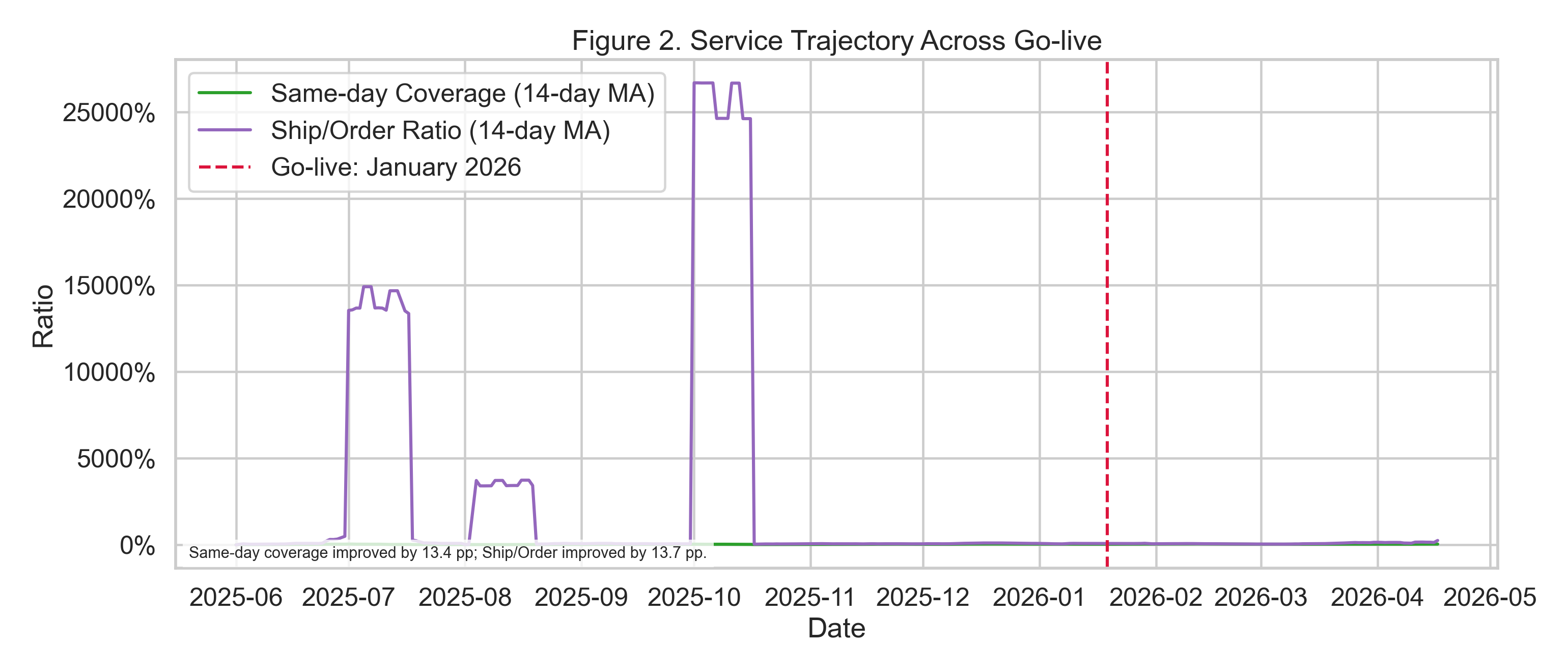}
\caption{Service trajectory before and after go-live (14-day moving averages). Green shows Same-Day Coverage; purple shows Ship/Order Ratio. Both trend upward after go-live.}
\label{fig:service-traj}
\end{figure}

Figure~\ref{fig:daily-volumes} and Figure~\ref{fig:service-traj} provide the temporal view of system behavior. The first figure focuses on absolute requested versus shipped batch volumes; the second figure focuses on normalized service quality trajectories.

\noindent\textbf{Interpretation (Fig.~\ref{fig:daily-volumes}).} This figure should be read as a throughput alignment chart: the gap between the blue and orange lines represents same-day unserved demand pressure. The visible narrowing of this gap after go-live indicates stronger operational response to incoming demand.

\noindent\textbf{Interpretation (Fig.~\ref{fig:service-traj}).} This figure isolates structural service change from raw volume fluctuations. The upward shift in both curves confirms that gains are not merely demand-mix artifacts: both bounded fulfillment and total shipment throughput improve together.

\subsection{Store-Level Capacity Compliance}
\begin{figure}[htbp]
\centering
\includegraphics[width=0.82\linewidth]{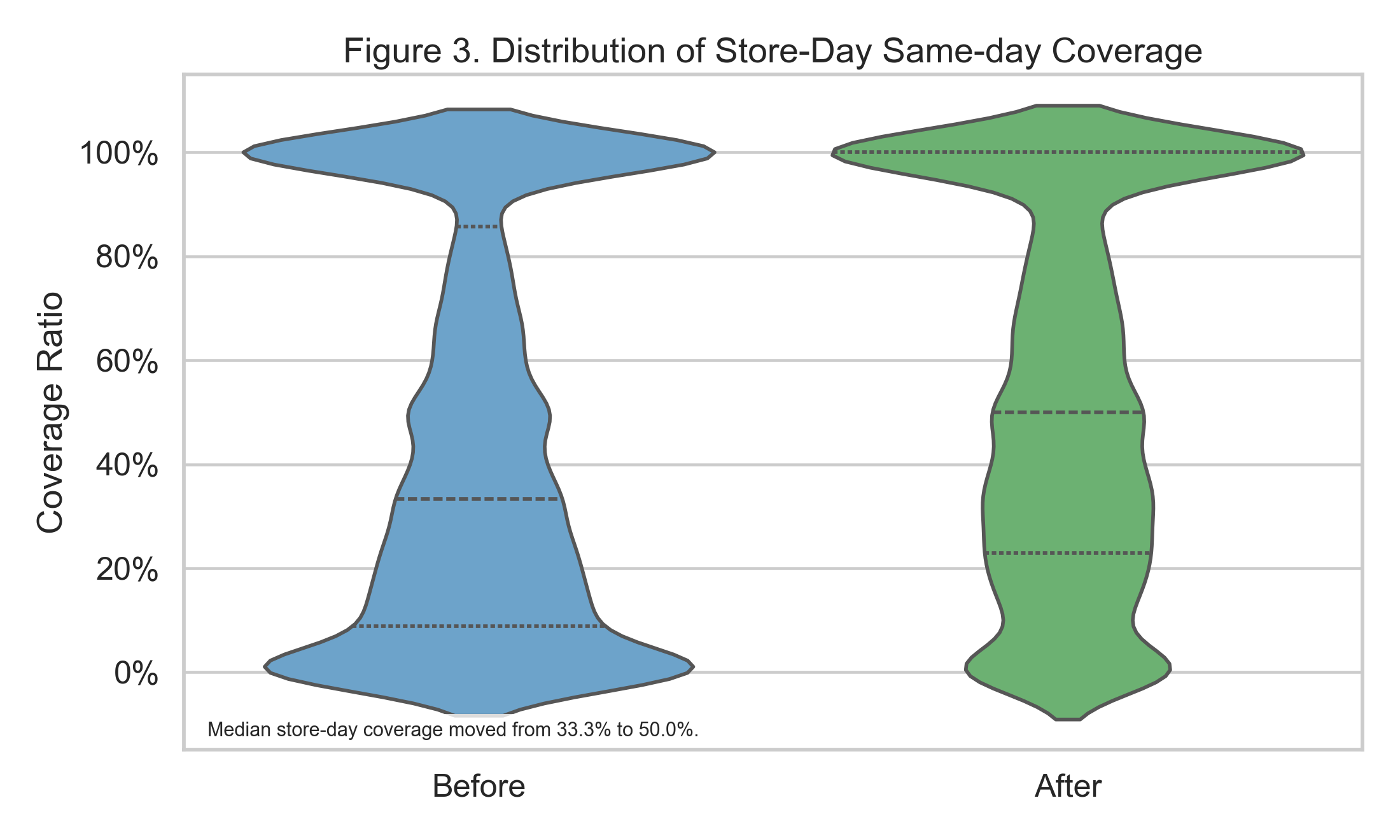}
\caption{Store-day coverage distribution before and after go-live. Quartile markers show that the center of the distribution shifts upward in the post period.}
\label{fig:coverage-dist}
\end{figure}

\begin{figure}[htbp]
\centering
\includegraphics[width=0.82\linewidth]{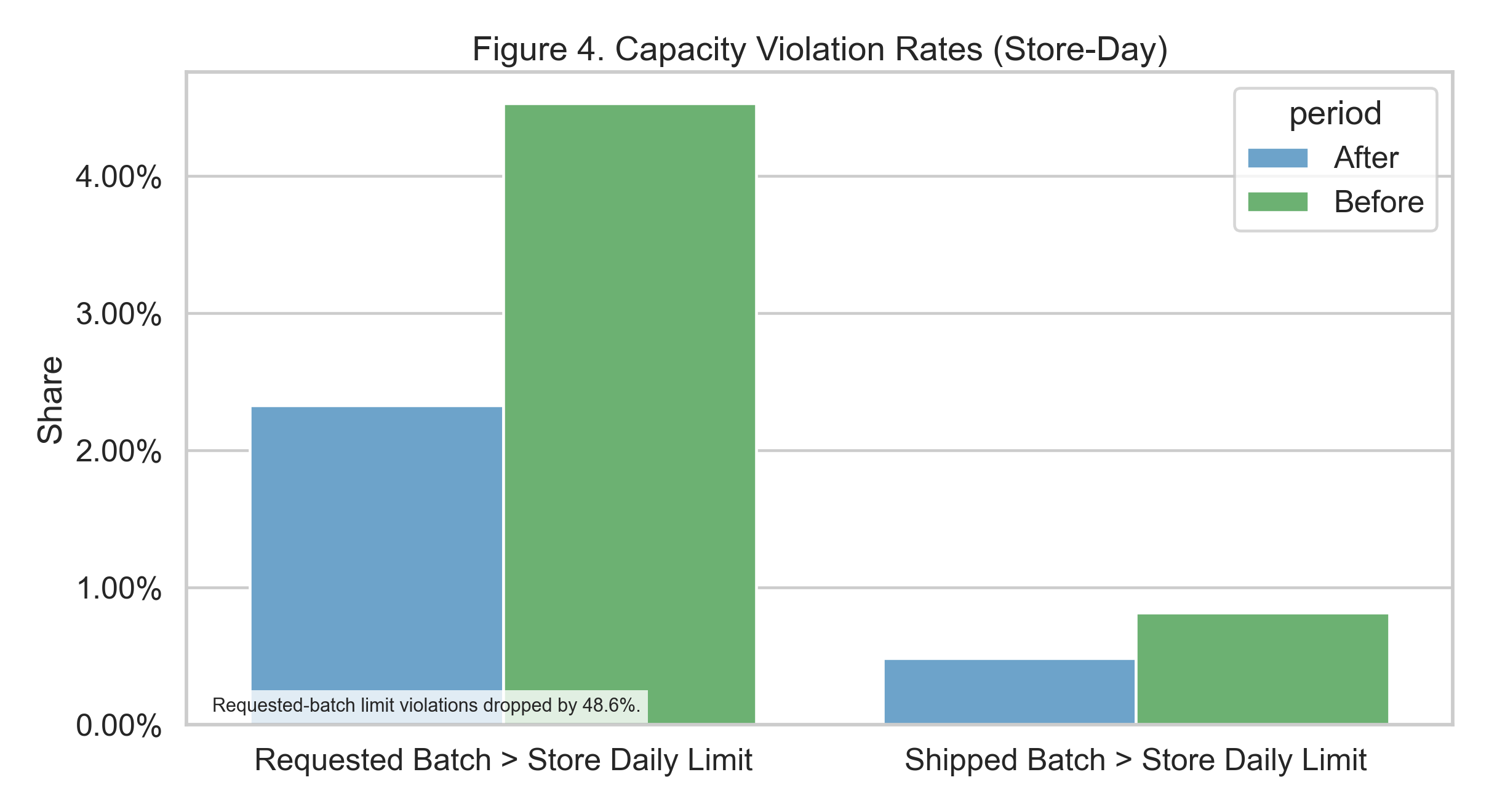}
\caption{Capacity violation shares at store-day level. Both requested-over-limit and shipped-over-limit frequencies decline after go-live.}
\label{fig:capacity-violations}
\end{figure}

Capacity compliance improved markedly after deployment. Figure~\ref{fig:coverage-dist} describes distributional movement of store-day coverage, while Figure~\ref{fig:capacity-violations} quantifies violation-rate reduction.

\noindent\textbf{Interpretation (Fig.~\ref{fig:coverage-dist}).} The post-period violin has higher central mass, indicating that coverage gains are broad-based across stores rather than concentrated in a few outliers.

\noindent\textbf{Interpretation (Fig.~\ref{fig:capacity-violations}).} The dominant reduction is on the requested-over-limit side, which is expected because the DSS applies limit-aware filtering before dispatch finalization. Shipped-over-limit also decreases, indicating improved execution discipline.

\subsection{Warehouse and 3PL Decomposition}
\begin{figure}[htbp]
\centering
\includegraphics[width=0.78\linewidth]{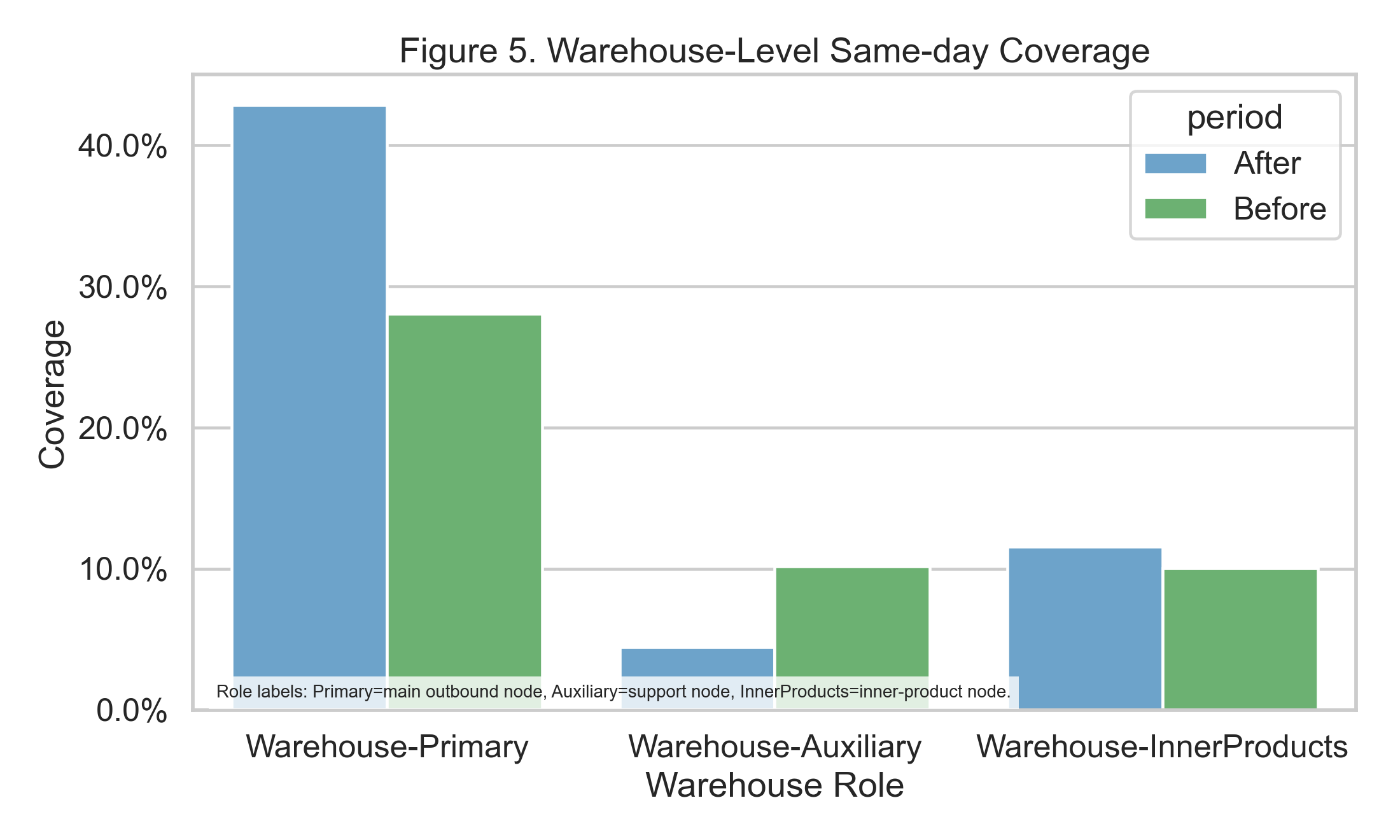}
\caption{Warehouse-level same-day coverage before and after go-live using role labels. Warehouse-Primary denotes the main outbound node, Warehouse-Auxiliary the support node, and Warehouse-InnerProducts the inner-product node.}
\label{fig:warehouse-cov}
\end{figure}

\begin{figure}[htbp]
\centering
\includegraphics[width=0.98\linewidth]{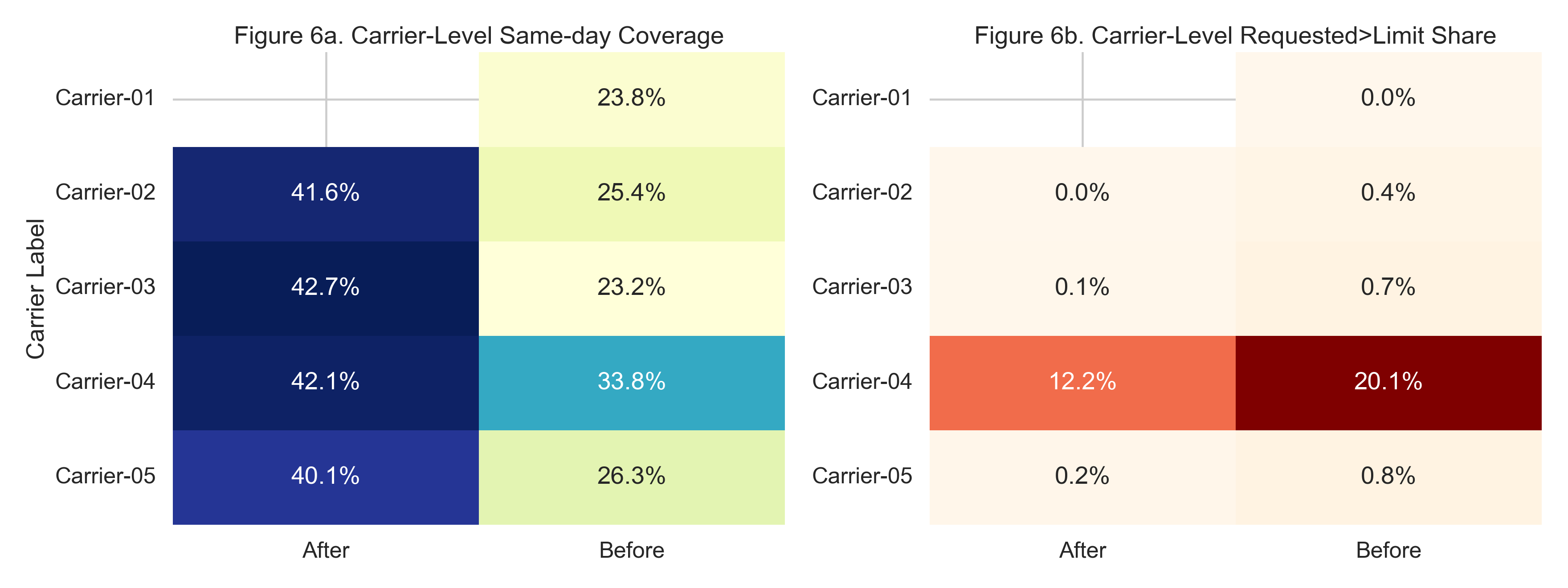}
\caption{Carrier decomposition: same-day coverage (left) and requested-over-limit share (right). Carriers are displayed as Carrier-01, Carrier-02, etc.}
\label{fig:pl-heatmap}
\end{figure}

Figure~\ref{fig:warehouse-cov} and Figure~\ref{fig:pl-heatmap} decompose the total effect by source-warehouse roles and carrier groups.

\noindent\textbf{Interpretation (Fig.~\ref{fig:warehouse-cov}).} The largest absolute uplift occurs in Warehouse-Primary, which also carries the highest total volume share; this indicates that global gains are not driven only by low-volume edge nodes.

\noindent\textbf{Interpretation (Fig.~\ref{fig:pl-heatmap}).} The left panel should be read as ``higher is better'' (coverage), while the right panel is ``lower is better'' (limit exceedance). Joint improvement across both panels suggests stable policy transfer across carrier lanes.

\subsection{Store Heterogeneity and Improvement Frontier}
\begin{figure}[htbp]
\centering
\includegraphics[width=0.92\linewidth]{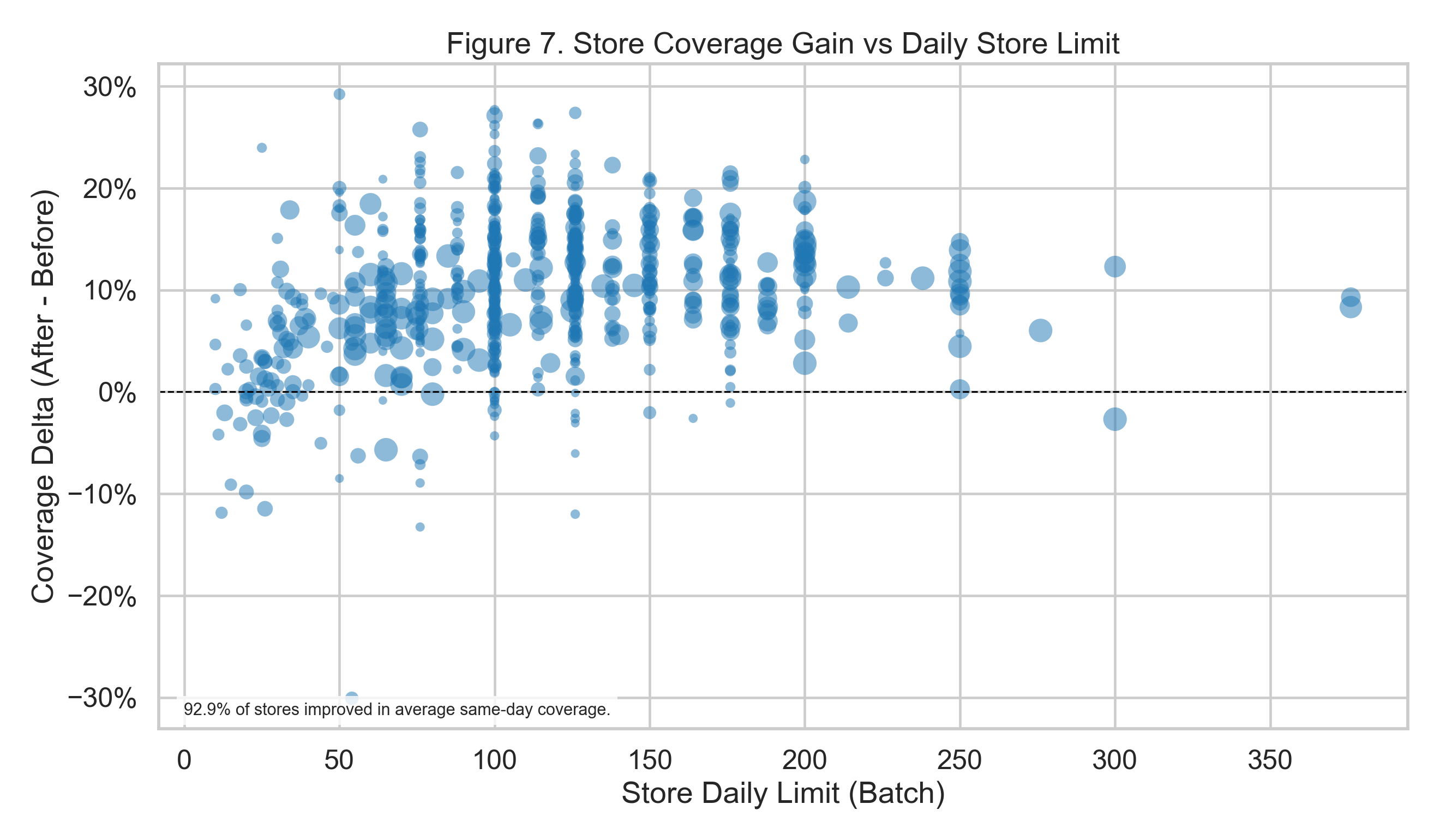}
\caption{Store-level coverage improvement versus daily store limit (bubble size proportional to total store volume).}
\label{fig:limit-vs-gain}
\end{figure}

\begin{figure}[htbp]
\centering
\includegraphics[width=0.82\linewidth]{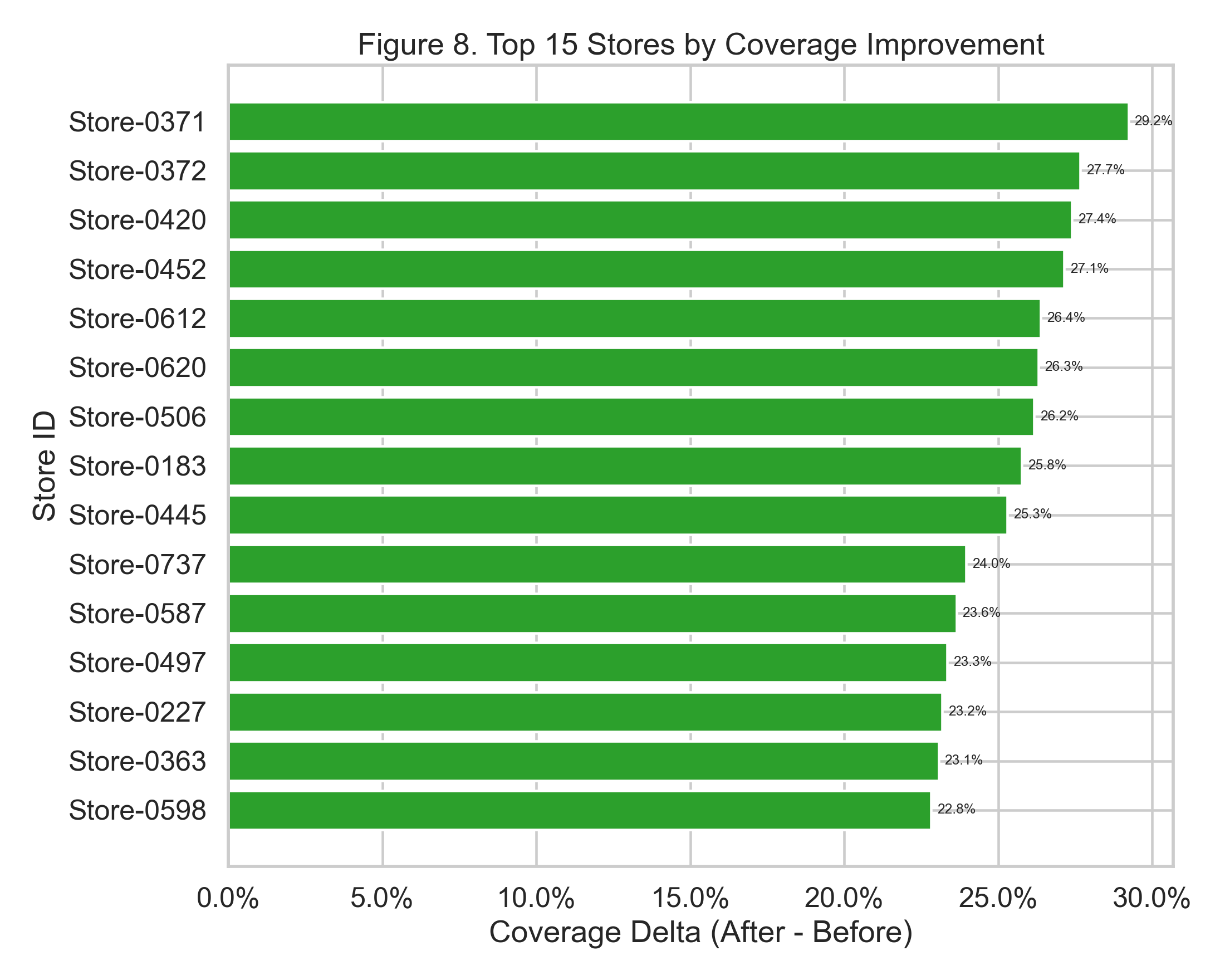}
\caption{Top 15 stores by same-day coverage improvement (After - Before), shown with Store-XXXX labels.}
\label{fig:top15}
\end{figure}

Coverage gains are not uniform across stores. Figure~\ref{fig:limit-vs-gain} links improvement magnitude to store limit levels, and Figure~\ref{fig:top15} highlights the leading improvement frontier under Store-XXXX labels.

\noindent\textbf{Interpretation (Fig.~\ref{fig:limit-vs-gain}).} Improvements are observed across a wide limit range, indicating that the DSS does not benefit only one store-capacity segment. The large bubbles above zero show that high-impact stores also improved.

\noindent\textbf{Interpretation (Fig.~\ref{fig:top15}).} This chart exposes the improvement concentration tail. It is useful for prioritizing post-implementation audits and extracting store-level best practices from high-gain nodes.

\section{Discussion and Managerial Implications}
\label{sec:discussion}

The results show that uncertainty-aware daily allocation control can deliver meaningful service and compliance gains without relying on computationally heavy exact optimization. In contrast to studies that focus primarily on game-theoretic equilibrium or stochastic inventory coupling \cite{Bahuguna2026,Roushdy2026}, our contribution emphasizes operational deployment speed in a live retail setting while preserving multi-constraint feasibility. From a theoretical standpoint, the paper contributes a unified constraint representation together with explicit feasibility and complexity properties for a deployable heuristic class.

From a managerial perspective, three implications are central. First, planners gain a robust daily mechanism for absorbing volatility in orders and residual capacities through rapid reruns. Second, warehouse activation and priority controls improve governance by making sourcing precedence explicit and auditable. Third, cumulative limit filtering lowers capacity violations and reduces execution friction for downstream logistics teams.

This study also has limitations. The benchmark is an operational before-after baseline rather than a full algorithmic comparison against exact and meta-heuristic solvers on matched synthetic instances. Future research can extend the framework with scenario-based stress testing, robust/stochastic optimization layers, and sampled-instance benchmarking to deepen theoretical comparison while preserving real-time implementability.

\section{Conclusion}
\label{sec:conclusion}

In modern supply chain architectures, balancing deterministic execution times with robust payload optimization is critical. The proposed data-driven real-time scalable heuristic maps retail replenishment into a constrained MKP pipeline that remains tractable in real time. The before-after evaluation around the January 2026 go-live shows consistent operational gains: weighted Ship/Order ratio improves from 54.1\% to 67.8\%, weighted same-day coverage improves from 24.3\% to 37.8\%, average daily unserved batch decreases by 22.8\%, and order-over-limit store-days are reduced by 48.6\%. These outcomes indicate that integrating outbound warehouse selection, warehouse-priority ranking, and cumulative capacity filtering within a computationally efficient DSS framework provides both scalability and stronger execution quality in high-variance omnichannel environments. Future work will include controlled benchmarking against exact and meta-heuristic alternatives on matched test instances. Future extensions will evaluate predictive machine-learning models to update both business priority scalar ($w_i$) and warehouse precedence policies before allocation.

\section*{Conflict of Interest}
The author declares that there is no conflict of interest regarding the publication of this manuscript.

\bibliographystyle{elsarticle-num}

\end{document}